\documentclass{amsart}
\usepackage{amsfonts,amscd}

\newtheorem{theorem}{Theorem}[section]
\newtheorem{lemma}[theorem]{Lemma}
\newtheorem{corollary}[theorem]{Corollary}
\newtheorem{proposition}[theorem]{Proposition}
\theoremstyle{remark}

\theoremstyle{definition}
\newtheorem{definition}[theorem]{Definition}

\numberwithin{equation}{section} \makeatother

\DeclareMathOperator{\Kdb}{{\mathbb K}}

\DeclareMathOperator{\Ndb}{{\mathbb N}}

\begin{document}
\title[Hereditary subalgebras of operator algebras]{Hereditary subalgebras of
operator algebras}

\author{David P. Blecher, Damon M. Hay, {\protect \and} Matthew Neal}
\address{DAVID P. BLECHER, Department of Mathematics, University of Houston, Houston, TX
77204-3008}
\email{dblecher@math.uh.edu}
\address{DAMON M. HAY, Department of Mathematics, University of Houston, Houston, TX
77204-3008}
\email{dhay@math.uh.edu}
\address{MATTHEW NEAL, Math.\ and Computer Science Department, Denison University, Granville, OH 43023}
\email{nealm@denison.edu}
\thanks{*Blecher and Hay were partially supported by grant DMS 0400731 from
the National Science Foundation. Neal was supported 
by Denison University.}

\begin{abstract}
 In recent work of the second author,
a technical result was proved establishing a bijective
correspondence between certain open projections in a
$C^*$-algebra containing an operator algebra $A$,
and certain one-sided ideals of $A$.   Here we give several
remarkable consequences of this result.   These include
a generalization of the theory of hereditary subalgebras of
a $C^*$-algebra, 
and the solution of a ten year old
problem concerning the Morita equivalence of operator algebras.
In particular, the latter gives a very clean
generalization of the notion of Hilbert $C^*$-modules to
nonselfadjoint  algebras.
We show that an `ideal'  of a general operator space $X$  is
the intersection of $X$ with an `ideal' in any containing 
$C^*$-algebra or $C^*$-module.  Finally, we discuss 
 the noncommutative variant of the classical theory of 
`peak sets'.  \end{abstract}

\thanks{MSC (2000) 
Primary 46H10, 46L07, 46L85, 47L30, Secondary 32T40, 46A55, 46L08, 46L30, 47L50.}

\maketitle

\section{INTRODUCTION}   In \cite{HTh,H},
a technical  result was proved establishing a bijective 
correspondence between certain open projections in a
$C^*$-algebra containing an operator algebra $A$,
and certain one-sided ideals of $A$.   Here we give several 
remarkable consequences of this result.   These include 
a generalization of the theory of hereditary subalgebras of
a $C^*$-algebra, 
and the solution of a ten year old
problem concerning the Morita equivalence of operator algebras.
In particular, the latter yields the conceptually
cleanest generalization of
the notion of Hilbert $C^*$-modules to
nonselfadjoint  algebras.
We show that an `ideal'  of a general operator space $X$  is
the intersection of $X$ with an `ideal' in any containing
$C^*$-algebra or $C^*$-module.  Finally, we discuss
 the noncommutative variant of the classical theory of
`peak sets'.   If $A$ is a function algebra on a compact 
space $X$, then a {\em $p$-set} may be characterized
as a closed subset $E$ of $X$ such that for any open set $U$
containing $E$ there is a function in Ball$(A)$ which is $1$ on $E$,
and $< \epsilon$ in modulus outside of $U$.  We prove 
a noncommutative version of this result.      

  An {\em operator algebra} is a closed algebra of operators on
a Hilbert space; or equivalently a closed subalgebra of a $C^*$-algebra.
We refer the reader to \cite{BLM} for the basic theory of 
operator algebras which we shall need.   We say that
an operator algebra $A$ is {\em unital}
if it has an identity of norm $1$, 
and {\em approximately unital} if it 
has a contractive approximate identity (cai).  
A {\em unital-subalgebra} of a $C^*$-algebra $B$
 is a closed subalgebra containing $1_B$.
In this paper we will often work 
with closed right ideals $J$ of an operator algebra $A$ 
possessing a contractive left approximate identity (or left cai)
for $J$.   For brevity we will call these 
{\em r-ideals}.
The matching class of left ideals
with right cai will be called {\em $\ell$-ideals}, but these will
not need to be mentioned much for reasons of symmetry.  
In fact  r-ideals are exactly the {\em right $M$-ideals} 
of $A$ if $A$ is approximately unital \cite{BEZ}.  For 
$C^*$-algebras r-ideals are precisely the right 
ideals, and there is an obvious bijective correspondence
between r-ideals and $\ell$-ideals, namely $J \mapsto J^*$.
For nonselfadjoint operator algebras it is not at all clear
that there is a bijective correspondence between r-ideals and $\ell$-ideals.
In fact there is, but this seems at present
to be a deep result, as we shall see.
It is easy to see that there is a bijective correspondence between r-ideals 
$J$ and certain projections $p$ in the second dual $A^{**}$ 
(we recall that $A^{**}$ is also an operator algebra 
\cite[Section 2.5]{BLM}).  This 
bijection takes $J$ to its {\em left support 
projection}, namely the weak* limit of a left cai
for $J$; and conversely takes $p$ to the right ideal
$p A^{**} \cap A$.   The main theorem from \cite{H},
which for brevity we will refer to as {\em Hay's theorem}, states
that if $A$ is a unital-subalgebra of a $C^*$-algebra $B$
then  the projections $p$ here may be characterized as the 
projections in $A^{\perp \perp}$ which are open in $B^{**}$
in the sense of e.g.\ \cite{Ake,Ped}.   
Although this result sounds innocuous, 
its proof is presently quite technical and lengthy, and 
uses the noncommutative Urysohn lemma 
\cite{Ake2} and various nonselfadjoint analogues of it. 
  One advantage of this condition is that it has a left/right 
symmetry, and thus it leads naturally into a theory of 
hereditary subalgebras ({\em HSA's} for short)
of general operator algebras.   For commutative 
$C^*$-algebras of course HSA's are precisely the
closed two-sided ideals.  For 
noncommutative $C^*$-algebras the hereditary subalgebras are
the intersections
of a right ideal with its canonically associated left ideal
\cite{Ef,Ped}.
They are also the selfadjoint `inner ideals'.   (In this paper, we say that a
subspace $J$ of an algebra $A$ is an
inner ideal if $J A J \subset J$.
Inner ideals in this sense are sometimes called `hereditary subalgebras' in
the literature, but we will reserve the latter term for 
something more specific.)   The fact that 
HSA's of $C^*$-algebras are the selfadjoint inner ideals
follows quickly
from Proposition \ref{her15} below and its proof, or it 
can be deduced from \cite{ER1}.
HSA's play some of the role that two-sided ideals play 
in the commutative theory.  Also, their usefulness stems in large part
because many important properties of the algebra pass to hereditary 
subalgebras
(for example, primeness or primitivity).

We now summarize the content of our paper.  In Section 2, we use 
Hay's theorem to generalize some of the 
$C^*$-algebraic theory of HSA's.  Also
in Section 2 we use our results to give a solution\footnote{
An earlier attempt to solve this problem was made in \cite{K}.}
 to a problem raised in
\cite{B}.    In that paper an operator algebra $A$ was said to have
{\em property  (${\mathcal L}$)} if it has a
left cai $(e_t)$ such that $e_s e_t \to e_s$ with $t$ for each $s$.
It was asked if every operator algebra with a
left cai has property (${\mathcal L}$).   
As an application of this, in Section 3 we settle
a problem going back to the early days of the theory of strong
Morita equivalence 
of nonselfadjoint
operator algebras.  This gives a very clean generalization of
the notion of Hilbert $C^*$-module to such algebras.
In Section 4, we generalize to
nonselfadjoint algebras the connections between HSA's, 
weak* closed faces of the state space, and lowersemicontinuity.
We remark that facial structure in the algebra itself
has been looked at in the nonselfadjoint
literature, for example in
\cite{Kat} and references therein.
In Section 5 we show that every right $M$-ideal in any
 operator space $X$ is an 
intersection of $X$ with a canonical 
right submodule of any $C^*$-module (or `TRO') containing $X$.
Similar results hold for two-sided, or `quasi-', $M$-ideals.
This generalizes to arbitrary 
operator spaces the theme from e.g.\ Theorem \ref{cher1} below, and 
from \cite{H}, that r-ideals (resp.\ 
HSA's) are very tightly related to 
matching right ideals (resp.\
HSA's) in a containing $C^*$-algebra.
In the final Section 6 we discuss connections
with the peak and $p$-projections
introduced in \cite{HTh,H}.   The motivation
for looking at these objects is to attempt to generalize
the tools of peak sets and `peak interpolation' from 
the classical theory of function algebras (due to 
Bishop, Glicksberg, Gamelin, and others).   In particular,
we reduce the main open question posed
in \cite{H},
namely whether the $p$-projections coincide with
the support projections of r-ideals,
to a simple sounding question about 
approximate identities: If $A$ is an approximately unital operator algebra
then does $A$ have an approximate identity
of form $(1- x_t)$ with $x_t \in {\rm Ball}(A^1)$?
Here $1$ is the identity of the unitization
$A^1$ of $A$.   
 We imagine that the answer to this is in the negative. 
 We also show that $p$-projections
are exactly the closed projections
satisfying 
the `nonselfadjoint Urysohn lemma' or `peaking'
property discussed at the beginning
of this introduction.  Thus even if the question above
turns out in the negative, these projections should play an important
role in future `nonselfadjoint interpolation theory'.      
 
Hereditary subalgebras of not necessarily selfadjoint unital
operator algebras have previously been considered in the
papers \cite{MZ,ZZ} on inner ideals.   We thank Lunchuan Zhang for
sending us a copy of these papers.
Another work that has a point of contact
with our paper is the unpublished note \cite{Kun}.
Here {\em quasi-$M$-ideals},
an interesting variant of the one-sided $M$-ideals of  Blecher, Effros,
and Zarikian \cite{BEZ}
were defined.
Kaneda showed that the product $R L$
of an r-ideal and an $\ell$-ideal in an approximately unital operator
algebra A is an inner ideal (inner ideals are called `quasi-ideals' there),
and is a quasi-$M$-ideal.  It is also  noted there that
in a $C^*$-algebra $A$, the following three are the same:
quasi-$M$-ideals in
A, products $R L$ of an r-ideal and an $\ell$-ideal, and
inner ideals (see also \cite{ER1},
particularly Corollary 2.6 there).
Hereditary subalgebras
in the sense of our paper were not considered in \cite{Kun}.
We thank Kaneda for permission to
describe his work here and in Section 5.   
  
Some notations: In this paper, all projections are orthogonal 
projections.  If $X$ and $Y$ are sets (in an operator algebra say)
then we write
$X Y$ for the {\em norm closure} of the span of terms
of the form $x y$, for $x \in X, y \in Y$.  
The second dual $A^{**}$ of an operator algebra $A$
is again an operator algebra, and the first dual
$A^*$ is a bimodule over $A^{**}$ via the 
actions described, for example, on the bottom of 
p.\ 78 of \cite{BLM}.       
A projection $p$ in the second dual of a $C^*$-algebra $B$
is called {\em open} it is the sup of an increasing net of 
positive elements of  $B$.
Such projections $p$ are in a bijective correspondence with 
the right ideals $J$ of $B$, or with
the HSA's (see \cite{Ped}).   It is well known, and easy
to see, that $p$ is open iff there is a net $(x_t)$ in $B$
with $x_t \to p$ weak*, and $p x_t = x_t$.    
We recall that
TRO's are essentially the same thing as Hilbert
$C^*$-modules, and may be viewed as
closed subspaces $Z$ of $C^*$-algebras with the
property that $Z Z^* Z \subset Z$.  See e.g.\ \cite[Section 8.3]{BLM}.
Every operator space $X$ has a `noncommutative 
Shilov boundary' or `ternary envelope' 
$(Z,j)$ consisting of a TRO $Z$ and a
complete isometry $j : X \to Z$ whose 
range `generates' $Z$.   This ternary envelope has a
universal property which may be found in \cite{Ham,BLM}:
For any complete isometry $i : X \to Y$ into a
TRO $Y$, whose
range `generates' $Y$, there exists a (necessarily unique and
surjective) `ternary morphism' $\theta : Y \to Z$ such that
$\theta \circ i = j$.
 If $A$ is
an approximately unital operator algebra then
the noncommutative Shilov boundary is written as $C^*_e(A)$ (see e.g.\ 
\cite[Section 4.3]{BLM}), and was first introduced
by Arveson \cite{SOC}.

\section{HEREDITARY SUBALGEBRAS}

Throughout this section $A$ is an operator algebra (possibly not 
approximately unital).  Then $A^{**}$ is an operator algebra.
   We shall say that
a projection  $p$ in $A^{**}$ is {\em
open in $A^{**}$} if $p \in (p A^{**} p  \cap A)^{\perp \perp}$.
In this case we also say that $p^{\perp}$ is
{\em closed in  $A^{**}$}, or is an {\em approximate
$p$-projection} (this notation was used in \cite{H}
since these projections have properties analoguous
to the {\em $p$-sets} in the theory of uniform algebras;
see e.g.\ \cite[Theorem 5.12]{H}).
Clearly these notions
are independent of
any particular $C^*$-algebra containing $A$.
If $A$ is a $C^*$-algebra then
these concepts coincide with the usual notion
of  open and closed projections (see e.g.\ \cite{Ake,Ped}).

\medskip

{\bf Example.} 
  Any projection 
$p$ in the multiplier algebra $M(A) \subset A^{**}$
is open in $A^{**}$, 
if $A$ is approximately unital.
Indeed $p A^{**} p \cap A = p A p$, and
if $(e_t)$ is a cai for $A$, then $p e_t p \to p$ weak*.

\medskip

If $p$ is open in $A^{**}$ then clearly $D = p A^{**} p \cap A$
is a closed subalgebra of $A$, and it has 
a cai by \cite[Proposition 2.5.8]{BLM}.  We call such a subalgebra
a  {\em hereditary subalgebra} of $A$ (or for brevity, a {\em HSA}). 
Perhaps more properly (in view of the next result) we should call these
`approximately unital 
HSA's', but for convenience we use the shorter term.  We say that
$p$ is the {\em support projection} of the HSA $D$; and it follows
by routine arguments that 
$p$ is the weak* limit of any cai from $D$.

\begin{proposition}[\cite{MZ,ZZ}]  \label{her15}
A subspace of an operator algebra $A$ is
a HSA iff it is an approximately unital inner ideal.
\end{proposition}  

\begin{proof}   We have already said that HSA's
are  approximately unital, and clearly they
are  inner ideals.  

If $J$ is an approximately unital inner ideal
then by \cite[Proposition 2.5.8]{BLM} we have 
that $J^{\perp \perp}$ is an algebra with
identity $e$ say.  Clearly 
$J^{\perp \perp} \subset e A^{**} e$.
Conversely, by a routine weak* density argument
$J^{\perp \perp}$ is an inner ideal, 
and so $J^{\perp \perp} = e A^{**} e$. 
 Thus $J = e A^{**} e \cap A$, and 
$e$ is open.  \end{proof}

We can often assume that the containing algebra
$A$ above is unital, simply
by adjoining a unit to $A$ (see \cite[Section 2.1]{BLM}).
 Indeed it follows from the last 
proposition that a subalgebra
$D$ of $A$ will be hereditary in the unitization $A^1$ 
iff it is hereditary in $A$.   

The following is a second (of many) characterization of
HSA's.   We leave the proof to the reader.

\begin{corollary}  \label{her8}  Let $A$ be an operator algebra
and suppose that $(e_t)$ is a net in ${\rm Ball}(A)$ such that
$e_t e_s \to e_s$ and $e_s e_t  \to e_s$ with $t$.
Then $\{ x \in A : x e_t  \to x , e_t x \to x \}$ is
a HSA of $A$.  Conversely,
every HSA of $A$ arises in this
way.
\end{corollary}

Note that this implies that
any approximately unital subalgebra $D$ of $A$
is contained in a HSA.

We next refine Hay's theorem from \cite{H}.

\begin{theorem}  \label{TGBTG}
Suppose that $A$ is an operator algebra (possibly not
approximately unital), 
and that $p$ is a projection in $A^{**}$.
The following are equivalent:
\begin{itemize}
\item [(i)]  $p$ is open in $A^{**}$.
\item [(ii)] $p$ is open as a projection in $B^{**}$,
if $B$ is a $C^*$-algebra containing $A$ as a subalgebra.
\item [(iii)] $p$ is the left support projection of an r-ideal
of $A$ (or, equivalently, $p$ is contained 
in $(p A^{**}  \cap A)^{\perp \perp}$).
\item [(iv)]   $p$ is the right support projection of an $\ell$-ideal of $A$.
\item [(v)] $p$ is the  support projection of a hereditary subalgebra
of $A$. 
  \end{itemize}
\end{theorem}

\begin{proof}  That (v) is equivalent to (i) is  just the  
definition of being open in $A^{**}$. 
Also (i) implies (ii) by 
facts about open projections mentioned in the introduction.
Supposing (ii), consider $A^1$ as a unital-subalgebra of $B^1$.
Then $p$ is open as a projection in $(B^1)^{**}$.
Since $p \in A^{\perp \perp}$ it follows from Hay's theorem
that $J = p (A^1)^{**} \cap A^1$ is an r-ideal
of $A^1$ with left support projection $p$.  If $x \in J$ then 
$x = p x \in (A^{\perp \perp} A^1) \cap A^1 \subset A^{\perp \perp}
 \cap A^1 = A$.   
Thus $J = p A^{**} \cap A$, and we have proved (iii).
Thus to complete the proof
it suffices to show that (iii) implies (i) (the equivalence with (iv)
following by symmetry).

(iii) $\Rightarrow$ (i) \
First assume that $A$ is unital, in which case
(iii) is equivalent to (ii) by Hay's theorem. 
We work in $A^*$.  As stated in the intoduction,  
$A^*$ is a right $A^{**}$-module via the 
action $(\psi \eta)(a) = \langle \eta a , \psi \rangle$
for $\psi \in A^*, \eta \in A^{**}, a \in A$.
Similarly it is a left $A^{**}$ module.  Let $q = p^\perp$,
a closed projection in $B^{**}$ for any $C^*$-algebra $B$
 generated by $A$.  We first claim that $A^* q = 
J^\perp$, where $J$ is the right ideal of $A$ corresponding to $p$,
and so $A^* q$ is weak* closed.   
  To see that $A^* q =
J^\perp$, note that clearly $A^* q\subset J^\perp$, since $J = p A^{**} \cap A$.
Thus if  $\psi \in J^\perp$, then $\psi q \in J^\perp$,
and so $\psi p \in J^\perp$ since $\psi = \psi p + \psi q$. 
However, if $\psi p \in J^\perp = (p A^{**})_\perp$,
then $\psi p \in (A^{**})_\perp = \{ 0 \}$.   Thus 
$\psi  = \psi q \in A^* q$.

Similarly, using the equivalence with (ii) here,
we have that $q A^*$ is weak* closed.   Now 
$q A^* + A^* q$ is the kernel of the projection 
$\psi \to p \psi p$ on $A^*$, and hence it is
norm closed.  By \cite[Lemma I.1.14]{HWW}, $q A^* + A^* q$ is 
weak* closed.   Claim: $(q A^* + A^* q)^\perp
= p A^{**} p$.    Assuming  this claim, note that 
$(q A^* + A^* q)_\perp \subset p A^{**} p \cap A$;
and $p A^{**} p \cap A  \subset (q A^* + A^* q)^\perp$,
so that $p A^{**} p \cap A = (q A^* + A^* q)_\perp$.
Thus $(p A^{**} p \cap A)^{\perp \perp} = p A^{**} p$, and the
proof is complete.  

In order to prove the claim, first note that it is clear that
$p A^{**} p \subset (q A^* + A^* q)^\perp$.  On the 
other hand, if $\eta \in (q A^* + A^* q)^\perp$ then 
write $\eta = p \eta p +  p \eta q + q \eta p + q \eta q$.
Thus $p \eta q + q \eta p + q \eta q \in (q A^* + A^* q)^\perp$.
In particular, applying this element to a functional $q \psi \in q A^*$
gives $$0 = \langle p \eta q + q \eta q ,  q \psi \rangle  =
\langle p \eta q + q \eta q , \psi \rangle , \qquad \psi \in A^* .$$
Thus $p \eta q + q \eta q = 0$, and left multiplying by $p$ shows that
$p \eta q = q \eta q = 0$.
Similarly $q \eta p = 0$.  Thus $\eta \in p A^{**} p$. 

Now assume that $A$ is nonunital.  If $J$ is the 
r-ideal, then  $J$ is an r-ideal
in $A^1$.
Thus by the earlier part, $p \in (p (A^1)^{**} p \cap A^1)^{\perp \perp}$.
 If $(e_t)$ is the  cai for $p (A^1)^{**} p \cap A^1$, 
then $e_t \to p$ weak*.
Since $p (A^1)^{**} \cap A^1 = J$ we have $p (A^1)^{**} p \cap A^1
\subset J \subset A$.  Thus $e_t \in p A^{**} p \cap A$, 
and so $p \in (p A^{**} p \cap A)^{\perp \perp}$.
Note too that the above shows that $p (A^1)^{**} p \cap A^1 =
p A^{**} p \cap A$.      \end{proof}

{\bf Remarks.}  1)  It is clear from the above
that a sup of open projections in $A^{**}$ is open in $A^{**}$.
From this remark, it is easy to give an alternative proof of
a result from \cite{BZ} which states that the closure of the
span of a family of r-ideals, again is an
r-ideal.

2)  If $A$ is approximately unital then one can add to the 
characterization of open projections in the theorem,
the condition that $A^* p^\perp$ is weak* closed in $A^*$. 
The second paragraph of the proof above shows one 
direction of this.  Conversely, if $A^* p^\perp$ is weak* closed, 
then $A^* p^\perp = J^\perp$ for a subspace $J$ of 
$A$  such that $J^{\perp \perp} = (A^* p^\perp)^\perp = p A^{**}$.
Thus $p$ is the support projection of the r-ideal
$A \cap p A^{**} = J$.

3)  A modification of part of the proof of the theorem shows
that if $A$ is approximately unital and if
$p, r$ are open projections in $A^{**}$ 
then 
$(p A^{**} r \cap A)^{\perp \perp}
= p A^{**} r$.   Note that $p A^{**} r \cap A$ is an inner ideal of
$A$.  Such subspaces are precisely the intersection of an 
r-ideal and an $\ell$-ideal.  
 
\begin{corollary} \label{propL}  Every operator algebra with a
left cai has property (${\mathcal L}$).
\end{corollary} 

\begin{proof}  
Let $C$ be an  operator algebra with a
left cai, and let $A$ be its unitization.
Then $C$ is an r-ideal in $A$, and the left support 
projection $p$ of $C$ in $A^{**}$ is a weak* limit 
of the left cai.   Also, $C = p A^{**} \cap A$.
By Theorem \ref{TGBTG}, we have 
$p \in (p A^{**} p  \cap A)^{\perp \perp}$,
and
$p A^{**} p  \cap A$ is a closed subalgebra
of $C$ containing a cai $(x_t)$ with 
$x_t \to p$ weak*.   If $J = \{ a \in A
:  x_t a \to a \}$ then $J$ is a  right ideal of $A$
with support projection $p$, so that $J = C$.  
Hence $C$ has property (${\mathcal L}$).
\end{proof} 

Some implications of this result are mentioned in
\cite{B}, however our main application appears in the 
next section.

In the following, we use some notation introduced in \cite{B}.
Namely, if $J$ is an operator algebra with a
left cai $(e_t)$ such that $e_s e_t \to e_s$ with $t$,
then we set ${\mathcal L}(J) = \{ a \in J : a e_t \to a \}$.
This latter space does not depend on the particular $(e_t)$,
as is shown in \cite{B}.

\begin{corollary}  \label{her3}
A subalgebra of an operator algebra $A$
is hereditary if and only if
it equals ${\mathcal L}(J)$ for an r-ideal $J$ of
$A$.
Moreover the correspondence $J \mapsto {\mathcal L}(J)$
is a bijection from the set of r-ideals of $A$ onto
the set of HSA's of $A$.  The inverse of
this bijection is the map
$D \mapsto D A$.  Similar results hold
for the $\ell$-ideals of $A$.
\end{corollary}
 
\begin{proof}
If $D$ is a HSA of $A$ then 
by Corollary \ref{her8} we have 
$D = \{ x \in A : x e_t  \to x , e_t x \to x \}$,
and $(e_t)$ is the cai for $D$.  Set $J = \{ x \in A : e_t x \to x \}$,
an r-ideal with $D = {\mathcal L}(J)$.
 
Conversely, if $J$ is an r-ideal then 
by Corollary \ref{propL}, 
we can choose a left cai $(e_t)$ of $J$ with the
property that $e_s e_t \to e_s$ with $t$.
Then $D = \{ x \in A :  x e_t  \to x , e_t x \to x \}$
is an HSA by Corollary \ref{her8}, and 
$D = {\mathcal L}(J)$.  Note that ${\mathcal L}(J) A
\subset J$, and conversely if $x \in J$ then
$x = \lim_t \, e_t x \in {\mathcal L}(J) A$.
Thus $J = {\mathcal L}(J) A$.  This shows that
$J \mapsto {\mathcal L}(J)$ is one-to-one.  The last
paragraph shows that it is onto. 
 \end{proof}

\begin{corollary} \label{her7}  If $D$ is a
hereditary subalgebra of an operator
algebra $A$, and if $J = D A$ and $K = A D$,
then $J K = J \cap K = D$.
\end{corollary}
 
\begin{proof}   Clearly $J K \subset J \cap K$.
Conversely, if $x \in J \cap K$ and $(e_t)$ is the cai for $D$
then $x = \lim_t x e_t \in J K$.  So $J K = J \cap K$ 
(see also e.g.\ \cite[Proposition 6.2]{FSIO} and \cite[Lemma 1.4.1]{Kth}).
Clearly $J K \subset D$ since $D$ is an inner 
ideal.  Conversely, $D = D^4 \subset J K$.
   \end{proof}        

\begin{theorem}  \label{cher1}
If $A$ is a closed subalgebra of a $C^*$-algebra $B$ then there 
is a bijective correspondence between r-ideals of $A$ and
right ideals of $B$ with left support in $A^{\perp \perp}$.
Similarly, there
is a bijective correspondence between HSA's of $A$ and 
 HSA's of $B$ with support in $A^{\perp \perp}$.
The correspondence takes an r-ideal (resp.\ HSA) $J$ of 
$A$ to $JB$ (resp.\ $J B J^*$).   The inverse
bijection is simply intersecting with $A$.    
\end{theorem}

\begin{proof}  We leave the proof of this to the reader, using 
the ideas above 
(and, in particular, Hay's theorem).  At some
point an appeal to \cite[Lemma 2.1.6]{BLM} might be necessary.    
\end{proof}

In the $C^*$-algebra case the correspondence between r-ideals
and $\ell$-ideals has a simple formula: $J \mapsto J^*$. 
For nonselfadjoint algebras $A$, one formula setting
up the same correspondence is $J \mapsto A {\mathcal L}(J)$.
It is easy to see from the last theorem that,
for subalgebras $A$ of a $C^*$-algebra $B$, this 
correspondence becomes $J \mapsto BJ^* \cap A$.
Here $J$ is an r-ideal; and notice that $BJ^* \cap A$ also equals
$B D^* \cap A$, where $D$ is the associated HSA 
of $A$
(we remark that by \cite[Lemma 2.1.6]{BLM} it is easy to see 
that $B D^* = B D$).   This allows
us to give another description of ${\mathcal L}(J)$ as 
$J \cap BJ^*$.

\begin{theorem}  Suppose $D$ is a hereditary subalgebra of an 
approximately unital operator algebra $A$. Then every
$f \in D^*$ has a unique 
Hahn-Banach extension to a functional in $A^{\ast}$ (of the same norm).
\end{theorem}

\begin{proof}
Let $g$ and $h$ be two such extensions.  Since $D=pA^{\ast\ast}p \cap A$ 
for an open projection $p$, it
is easy to see that $pgp = php$.  Since
$\|g\|=\|pgp\|=\|php\|=\|h\|$, we need only show that $g=pgp$ and similarly $h=php$.
Consider $A^{\ast\ast}$ as a unital-subalgebra of a W*-algebra $B$.
Since the canonical projection from $B$ onto $pBp +(1-p)B(1-p)$ is contractive,
and since $\Vert p b p + (1-p) b (1-p) \Vert = \max \{
\Vert p b p \Vert , \Vert (1-p) b (1-p) \Vert \}$ 
for $b \in B$, it is easy to 
argue that 
 \[
\|g\| \geq \|pgp +(1-p)g(1-p)\| = \|pgp\|+\|(1-p)g(1-p)\|\geq \|g\| .
\]
Hence, $(1-p)g(1-p)=0$.  Since $g = pgp + pg(1-p)+pg(1-p) + (1-p)g(1-p)$,
it suffices to show that $pg(1-p)+pg(1-p)=0$.  To this end,
we follow the proof in
Proposition 1 of \cite{FR}, which proves the analogous result for JB*-triples.
For the readers convenience,
we will reproduce this pretty argument in 
our setting, adding a few more details.
Of course $B$ is a JB*-triple.   We will use the notation
$pBp=B_{2}(p)$, $pB(1-p)+(1-p)Bp =B_{1}(p)$, and $(1-p)B(1-p) =B_{0}(p)$. 
For this proof only, we will write $x^{2n+1}$ for $x(x^{\ast}x)^{n}$
(this unusual notation is used in the JB*-triple literature).
   In Lemma 1.4 of
\cite{FR}, it is proved that, for $x \in B_{2}(p) \cup B_{0}(p)$, 
$y \in B_{1}(p)$, and $t > 0$,
\begin{eqnarray}
\; \; (x+ty)^{3^{n}}
=x^{3^{n}}+t2^{n}D(x^{3^{n-1}},x^{3^{n-1}}) \cdots D(x^{3},x^{3})D(x,x)y+
O(t^{2}).  \label{J}
\end{eqnarray}
\noindent where, in our setting, $D(w,w)$ is the operator
$D(w,w)z=(ww^{\ast}z+z w^{\ast} w)/2$ on $B$.
 Here, $O(t^{2})$ denotes
a polynomial in $x$, $y$, and $t$,
with all terms at least quadratic in $t$.  This polynomial has 
a certain number of terms that depends only on $n$, and
the coefficients of the monomials in $x, y$ and $t$ also 
depend only on $n$.    

Choose $y \in B_{1}(p) \cap
A^{\ast\ast}$.   We may assume  that $\|g\|=1$, $g(y)\geq 0$, and $\|y\|
\leq 1$.  Given $\epsilon > 0$, we choose $x
\in D$ with $\|x\|=1$ and $f(x) \geq 1-\epsilon$.   Then, for $t > 0$, we
have 
\[
\|x+ty\| \geq g(x+ty)=g(x)+t g(y) \geq 1-\epsilon +tg(y).
\]
Thus, by (\ref{J}) above, and the fact that $\|x\| \leq 1$,
\begin{eqnarray*}
(1- \epsilon +tg(y))^{3^{n}} \leq \|x+ty\|^{3^{n}} &=& \|(x+ty)^{3^{n}}\|\\
&\leq & \|x^{3^{n}}\|+t2^{n}\|y\|+ \Vert O(t^{2}) \Vert \\
&\leq &  1 + t2^{n}\|y\|+ p(t) ,
\end{eqnarray*}
where $p(t)$ is a polynomial in $t$ with all terms at least degree 2,
and coefficients which depend only on $n$ and $\Vert y \Vert$.
Letting $\epsilon \rightarrow 0$ we have 
$(1 +tg(y))^{3^{n}} \leq 1 + t2^{n}\|y\|+ p(t)$, and so 
\[
1 + 3^{n}tg(y)  \leq 
1+t2^{n}\|y\|+ r(t), 
\]
where $r(t)$ is a polynomial with the same properties as $p$,
and in particular has all terms at least degree 2.
Dividing by $3^{n} t$, we obtain
\[
g(y) \leq
\left(\frac{2}{3}\right)^{n}\|y\|+\frac{r(t)}{t 3^{n}} .   
\]
Letting $t \rightarrow 0$ and then 
$n \rightarrow \infty$, we see that $g(y)=0$.  Hence
$pg(1-p)+(1-p)gp = 0$ as desired.
\end{proof}

One might hope to improve the previous theorem 
to address extensions of completely bounded maps from $D$ into $B(H)$.
Unfortunately, simple examples such as the one-dimensional
HSA $D$ in $\ell^\infty_2$ which is supported in the
first entry, with $f : D \to M_2$ taking $(1,0)$ to $E_{11}$,
shows that one needs to impose strong restrictions on the 
extensions.    This two dimensional example contradicts 
several theorems on unique completely
contractive extensions in the literature.
We found the following positive result after reading \cite{ZZ}.
Although some part of it
is somewhat tautological, it may be the best that one could hope for.
To explain the notation in (iii), if $A$ is an
approximately unital operator algebra and $B$ is a unital weak* closed
operator algebra,  
then we say that a bounded map $T : A \to B$ is {\em weakly nondegenerate}  if
the canonical weak* continuous extension $\tilde{T} : A^{**} \to B$ is unital.
By 1.4.8 in \cite{BLM} for example, this is equivalent to: 
$T(e_t) \to 1_B$ weak*  for some contractive approximate
identity $(e_t)$ of $A$; and is also equivalent to the same statement
with `some' replaced by `every'.

\begin{proposition} \label{unx} Let $D$ be an
 approximately unital subalgebra
of an approximately unital operator algebra $A$.
The following are equivalent:
\begin{itemize}
\item [(i)] $A$ is a hereditary subalgebra of $A$.
 \item [(ii)]  Every completely contractive unital 
map from $D^{**}$ into a unital operator algebra $B$, has
a unique completely contractive unital extension from
$A^{**}$ into $B$. 
 \item [(iii)]  Every completely contractive weakly  nondegenerate
map from $D$ into a unital weak* closed operator algebra $B$
has a unique completely contractive  weakly  nondegenerate extension from $A$
into $B$.
\end{itemize} 
 \end{proposition}

\begin{proof}  We are identifying $D^{**}$ with $D^{\perp \perp}
\subset A^{**}$.  Let $e$ be the identity of $D^{**}$.

(ii) $\Rightarrow$ (i)  \ If (ii)   
holds,
then the identity map on $D^{**}$ extends  to a
unital complete contraction $S : A^{**} \to D^{**} \subset e A^{**} e$.  
The map $x \mapsto e x e$ on $A^{**}$
is also a completely contractive unital extension
of the inclusion map $D^{**} \to e A^{**} e$.  
It follows from the hypothesis that these maps
coincide, and so $e A^{**} e = D^{**}$, which implies 
that $D$ is a HSA.   

(i) $\Rightarrow$ (ii)  \ If $D$ is a HSA, then extensions
of the desired kind exist by virtue of the canonical
projection from $A^{**}$ onto $D^{\perp \perp}$.  For the 
uniqueness, suppose that $\Phi$ is  
such an extension of a completely contractive unital
map $T : D^{\perp \perp} \to B$. 
Since  $e$ is
an orthogonal projection in $A^{**}$, it follows from
the last remark in 2.6.16 in \cite{BLM} that 
$$T(e x e) =  \Phi(e x e) = \Phi(e) \Phi(x) \Phi(e) = \Phi(x) , \qquad x \in A^{**}. $$
Hence (ii) holds.

Inspecting the proof above shows that (i) is
equivalent to the variant of (ii) where $B$ is
weak* closed and all maps are 
also weak* continuous.  Then the equivalence with (iii)
is easy to see using the facts immediately above the 
proposition statement, and also the bijective 
correspondence between complete contractions
$A \to B$ and weak* continuous complete contractions $A^{**} \to B$
(see 1.4.8 in \cite{BLM}).        \end{proof}

\section{APPLICATION: A GENERALIZATION OF $C^*$-MODULES}

In the early 1990's, the first author together with Muhly and Paulsen
generalized Rieffel's strong
Morita equivalence to nonselfadjoint operator algebras \cite{BMP}.
This  study was extended to include a
generalization of Hilbert $C^*$-modules to nonselfadjoint algebras,
which were called {\em rigged modules} in  \cite{Bghm},
and {\em (P)-modules} in \cite{BMP}.  See \cite[Section 11]{Bnat} for 
a survey.   There are very many
equivalent definitions of these objects in these papers.   
The main purpose
of this section is to settle a problem going back to the
early days of this theory.  This results in the conceptually
clearest definition of rigged modules; and also tidies
up one of the characterizations of strong Morita equivalence.
The key tool we will use is
the Corollary  \ref{propL} to our main theorem from Section 2.

Throughout this section, $A$ is an approximately unital operator algebra.
For a positive integer $n$ we write $C_n(A)$ for
the $n \times 1$ matrices with entries in $A$, which
may be thought of as the first column of the
operator algebra $M_n(A)$.
 In our earlier work mentioned above
$C_n(A)$ plays the role of the prototypical
right $A$-rigged module, out of which all others
may be built via `asymptotic factorizations' similar
to the kind considered next.

\begin{definition} \label{rig}
An operator space $Y$ which is
also a right $A$-module is {\em $A$-Hilbertian}
if there exists a net  of positive integers $n_\alpha$,
and completely contractive $A$-module maps
$\varphi_\alpha : Y \to C_{n_\alpha}(A)$ and
$\psi_\alpha :  C_{n_\alpha}(A) \to Y$,
such that $\psi_\alpha \,
\varphi_\alpha \to I_Y$ strongly on $Y$.
  \end{definition}

The name `$A$-Hilbertian' is due to Paulsen around 1992, who 
suggested that these modules should play an important 
role in the Morita theory.
A few years later
the question became 
whether they coincide with the
rigged modules/(P)-modules from \cite{Bghm,BMP}.
This question appears explicitly in \cite[Section 11]{Bnat} for example,
and was discussed also several times 
in \cite[Chapter 4]{BMP} in terms of the 
necessity of adding further conditions
to what we called the `approximate identity property'.
Assuming for simplicity that $A$ is unital, one of the many equivalent
definitions of rigged modules is that they are the
modules satisfying Definition \ref{rig},
but that in addition $\varphi_\beta \psi_\alpha \,
\varphi_\alpha \to  \varphi_\beta$ in norm for 
each fixed $\beta$ in the directed set.   We were 
not able to get the theory going without this extra
condition.   Thus 
the open question referred to above may be restated
as follows: can one always replace
the given nets in Definition \ref{rig} with
ones which satisfy this  additional condition?  
The first author proved this if $A$ is a $C^*$-algebra in \cite{Bna}; indeed 
$A$-Hilbertian modules coincide with $C^*$-modules if $A$ is
a $C^*$-algebra.  A simpler proof of this
result due to Kirchberg is included in \cite[Theorem 4.24]{BMP}.

Although the `asymptotic factorization' in the definition above
is clean, it can sometimes
be clumsy to work with, as is somewhat illustrated
by the proof of the next result.

\begin{proposition} \label{prol}   Let $Y$ be an operator space
and  right $A$-module, such that
there exists a net  of positive integers $n_\alpha$,
$A$-Hilbertian modules $Y_\alpha$, and completely contractive $A$-module maps
$\varphi_\alpha : Y \to Y_\alpha$ and
$\psi_\alpha :  Y_\alpha \to Y$, such that $\psi_\alpha \,
\varphi_\alpha \to I_Y$ strongly.  Then $Y$ is
$A$-Hilbertian.
\end{proposition}

\begin{proof}
We use a net reindexing argument based on \cite[Lemma 2.1]{Bghm}.
  Suppose that
$\sigma^\alpha_\beta : Y_\alpha \to Z^\alpha_{\beta}$
and $\tau^\alpha_\beta : Z^\alpha_{\beta} \to Y_\alpha$,
are the `asymptotic factorization' nets corresponding to $Y_\alpha$.
We define a new directed set
$\Gamma$ consisting of 4-tuples
$\gamma = (\alpha, \beta, V, \epsilon)$,
where $V$ is a finite subset of $Y$, $\epsilon > 0$,
and such that $$\Vert \psi_\alpha \tau^\alpha_\beta  \sigma^\alpha_\beta
\varphi_\alpha(y) - \psi_\alpha  \varphi_\alpha(y) \Vert < \epsilon ,
\qquad y \in V .$$
This is a directed set with ordering
$(\alpha, \beta, V, \epsilon) \leq (\alpha',\beta',V',\epsilon')$
iff $\alpha \leq  \alpha', V \subset V'$ and $\epsilon' \leq \epsilon$.
(We recall that directed sets for nets make no essential 
use of the `antisymmetry' condition for the ordering, and we follow many
authors in not requiring this.)
Define $\varphi^\gamma = \sigma^\alpha_\beta \, \varphi_\alpha$
and $\psi^\gamma =  \psi_\alpha \tau^\alpha_\beta$, if
$\gamma = (\alpha, \beta, V, \epsilon)$.  Given $y \in Y$
and $\epsilon > 0$, choose $\alpha_0$ such that
$\Vert \psi_\alpha  \, \varphi_\alpha (y) - y \Vert < \epsilon$ whenever
$\alpha \geq \alpha_0$.  Choose $\beta_0$ such that
$\gamma_0 = (\alpha_0, \beta_0, \{ y \} , \epsilon) \in \Gamma$.
If $\gamma \geq \gamma_0$ in $\Gamma$ then
$$\Vert \psi^\gamma \, \varphi^\gamma(y)  - y \Vert
\leq \Vert \psi_\alpha \tau^\alpha_\beta  \sigma^\alpha_\beta
\varphi_\alpha(y) - \psi_\alpha  \varphi_\alpha(y) \Vert
+ \Vert \psi_\alpha  \varphi_\alpha(y) - y  \Vert < \epsilon' + \epsilon
\leq 2 \epsilon .$$
Thus $\psi^\gamma \,
\varphi^\gamma(y) \to y$,  and so $Y$ is 
$A$-Hilbertian.
\end{proof}

{\bf Remark.}     If desired, the appearance of the integers $n_\alpha$ in Definition
\ref{rig} may be avoided by the following trick.  
Let $C_\infty(A)$ be the space of columns $[x_k]_{k \in \Ndb}$, with $x_k \in A$,
 such that $\sum_k \, x^*_k x_k$ converges in $A$.
It is easy to see that $C_\infty(A)$ is $A$-Hilbertian, and for any 
 $m  \in \Ndb$ there is an obvious factorization of the identity map
on $C_m(A)$ through $C_\infty(A)$.   
It follows from this, and from the last Proposition, that Definition
\ref{rig} will be unchanged if all occurrences
of $n_\alpha$ there are replaced by $\infty$.

\begin{theorem} \label{main}   An operator space $Y$ which is
also a right $A$-module is a
rigged $A$-module if and only if it is  $A$-Hilbertian.
This is also equivalent to $Y$ having 
the `approximate identity property' of
{\rm \cite[Definition 4.6]{BMP}}.     \end{theorem}

\begin{proof}  Suppose that $Y$ is an operator space and a
right $A$-module which is $A$-Hilbertian.
It is easy to see that $Y$ is
an
operator  $A$-module, since the $C_{n_{\alpha}}(A)$ are,
and since $$\Vert [y_{ij}] \Vert
= \sup_\alpha \Vert [ \varphi_\alpha(y_{ij}) ]
\Vert = \lim_\alpha \Vert [ \varphi_\alpha(y_{ij}) ] \Vert
, \qquad [y_{ij}] \in M_n(Y) .$$
If 
$(e_t)$ is a cai for $A$ then the triangle inequality easily 
yields that for any $\alpha$,
\begin{eqnarray*}
\Vert y - y e_t \Vert  &  = & \Vert y - \psi_\alpha \varphi_\alpha(y) 
+ \psi_\alpha(\varphi_\alpha(y) - \varphi_\alpha(y) e_t) 
+ (\psi_\alpha \varphi_\alpha(y) - y) e_t \Vert \\
& \leq &  2 
\Vert y - \psi_\alpha \varphi_\alpha(y)) \Vert 
+ \Vert \varphi_\alpha(y) - \varphi_\alpha(y) e_t \Vert ,
\end{eqnarray*}
from which the nondegeneracy is easily seen.
Next, we reduce to the unital case.
Let $B = A^1$, the unitization of $A$.  Note that
$A$ is $B$-Hilbertian: the maps
$A \to B$ and
$B \to A$ being respectively the
inclusion, and left multiplication by elements in the cai
$(e_t)$.
Tensoring these maps with the identity map on
$C_{n_{\alpha}}$, we see that
$C_{n_{\alpha}}(A)$ is $B$-Hilbertian.
By Proposition \ref{prol},
$Y$ is $B$-Hilbertian.
By \cite[Proposition 2.5]{Bghm} it is easy to see that
$Y$ satisfies (the right module variant of)
\cite[Definition 4.6 (ii)]{BMP}.
By the results following that definition
we have that 
$C = Y \otimes_{hB} CB_B(Y,B)$ is
a closed right ideal in $CB_B(Y)$ which has a
left cai.  By
\cite[Theorem 2.7]{Bghm} or \cite[Theorem 4.9]{BMP}
we know that $CB_B(Y)$ is a unital operator
algebra.

By Corollary \ref{propL},
$C$ possesses a left cai $(v_\beta)$ such that
$v_{\gamma} v_\beta \to v_{\gamma}$ with $\beta$ for
each $\gamma$.    Let $D = \{ a \in C : a  v_\beta \to a
\}$, which is an operator algebra with cai.
Since the
 uncompleted algebraic tensor product $J$ of $Y$ with
$CB_B(Y,B)$ is
a dense ideal in $C$, and since
$J v_{\gamma} \subset D$ for each $\gamma$, it is
easy to see by the triangle inequality
that $J \cap D$ is a
dense ideal in $D$.  Thus we can rechoose a cai $(u_\nu)$
for $D$ from this ideal, if necessary by using
\cite[Lemma 2.1]{Bghm}.  This
cai will be a left cai for $C$ (e.g.\ see proof
of
Corollary \ref{propL}).  This implies that
$Y$ satisfies (a), and hence also (b), of
\cite[Definition 4.12]{BMP}.   That is,
$Y$ is a (P)-module, or equivalently a
rigged module, over $B$.   It is known that 
this implies that $Y$ is $A$-rigged.  One way to see this  
is to observe that by an application
of  Cohen's
factorization theorem as in \cite[Lemma 8.5.2]{BLM}, we have
$B_B(Y,B) = B_A(Y,A)$.
It follows that $Y$ satisfies  \cite[Definition 4.12]{BMP}
as an $A$-module, and hence
$Y$ is an $A$-rigged module.
That every rigged module
is $A$-Hilbertian follows from \cite[Definition 3.1]{Bghm}.
The equivalence with the `approximate identity property'
is essentially contained in the above argument.   
\end{proof}

This theorem
impacts only a small portion
of \cite{Bghm}.  Namely, that paper may now be improved
by replacing Definition 3.1 there by the modules
in Definition \ref{rig} above; and by tidying up
some of the surrounding exposition.
One may also now give alternative constructions of, for example,
the interior tensor product of rigged modules, by following
the idea in \cite[Theorem 8.2.11]{BLM}.

Similarly,  one may now tidy up
one of the characterizations of strong Morita equivalence.
By the above, what was called the `approximate identity 
property' in \cite[Chapter 4]{BMP} implies that the module 
is a (P)-module, and so in Theorems 4.21 and 4.23 in 
\cite{BMP} one may
replace conditions from Definition 4.12 with those in 4.6.
That is, we have the following improved 
characterization of the strong Morita equivalence of \cite{BMP}.
(The reader needing further details is referred to that source.) 

\begin{theorem}  \label{morth}  If $Y$ is a right $A$-Hilbertian 
module with the `dual approximate identity property' of \cite[Definition 4.18]{BMP},
then $Y$ implements a strong Morita equivalence between $A$ 
and the algebra $\Kdb_A(Y)$ of so-called `compact' operators on $Y$.
Conversely, every strong Morita equivalence
arises in such a way.  \end{theorem}

\medskip

{\bf Remark.}    The `dual approximate identity property'
mentioned in the theorem may also be phrased in terms of 
`asymptotic factorization'
of $I_A$ through spaces of the form $C_m(Y)$---this is
mentioned in \cite[p.\ 416]{Bghm} with a mistake that is
discussed in \cite[Remark 4.20]{BMP}.

\medskip

We refer the reader to \cite{Bghm} for the theory of rigged modules.  
It is easy to see using Corollary \ref{her7}
or Theorem \ref{main},  
that any hereditary subalgebra $D$ of an approximately unital
operator algebra $A$ gives rise to a rigged module.
Indeed, if $J = D A$, then
$J$ is a right rigged $A$-module, the 
canonical dual rigged module $\tilde{J}$ is just the
matching 
$\ell$-ideal $A D$, and the
operator algebra $\Kdb_A(J)$ of `compact operators' on $J$ is just $D$
 completely isometrically isomorphically.
From the theory of rigged modules \cite{Bghm} we know for example that
any completely contractive representation of $A$ induces
a completely contractive representation of $D$, 
and vice versa.  More
generally, any left operator $A$-module will give rise
to a left operator $D$-module by left tensoring with $J$,
and vice versa by left tensoring with $\tilde{J}$.  Since 
$J \otimes_{hA} \tilde{J} = D$ it follows that there is
an `injective' (but not in general `surjective') functor from 
$D$-modules to $A$-modules.
 
If $A$ and $B$ are approximately unital
operator algebras which are strongly Morita
equivalent in the sense of \cite{BMP},
then $A$ and $B$ will clearly be hereditary subalgebras of the
`linking operator algebra' associated with the
Morita equivalence \cite{BMP}.
Unfortunately, unlike the $C^*$-algebra case,
not every HSA $D$ of an operator algebra $A$
need be strongly Morita equivalent
to $A D A$.   One would also need a condition similar to
that of \cite[Definition 5.10]{BMP}.  Assuming the presence
of such an extra condition, it follows that the representation
theory for the algebra $A$ is `the same' as the representation
theory of $D$; as is always the case if one has a
Morita equivalence.

\medskip

{\bf Example.}
 If $a \in {\rm Ball}(A)$, for an
operator algebra $A$, let $D$ be the
closure of $(1-a) A (1-a)$.  Then it follows from
the later Lemma \ref{wilp} that $D$ is a   hereditary subalgebra
of $A$.  The associated r-ideal is $J$, the
closure of $(1-a) A$.  The dual rigged module
 $\tilde{J}$ is equal to the
 closure of $A (1-a)$, and $\Kdb_A(J) \cong D.$
It is easy to check that even for examples of this kind,
the $C^*$-algebra $C^*(D)$ generated by $D$
need not be a hereditary $C^*$-subalgebra
of $C^*(A)$ or $C^*_e(A)$.  For example, take $A$
to be the subalgebra of $M_2(B(H))$ consisting of all
matrices whose 1-1 and 2-2 entries
are scalar multiples of $I_H$, and whose 2-1 entry is $0$.
Let $a = 0_H \oplus I_H$.
In this case $D = (1-a) A (1-a)$ is one dimensional, and it
is not a HSA of $C^*(A)$.  Also $D$ is
not strongly Morita equivalent to $A D A$.
                    
\section{CLOSED FACES AND LOWERSEMICONTINUITY} 

Suppose that $A$ is an approximately unital operator algebra.
The state space
$S(A)$ is the set of functionals $\varphi \in A^*$  such that
$\Vert \varphi \Vert = \lim_t \varphi(e_t) = 1$,
if $(e_t)$ is a cai for $A$.  These are
all restrictions to $A$ of states on any $C^*$-algebra generated
by $A$.                                       
If $p$ is
a projection in $A^{**}$, then any $\varphi \in S(A)$ may be 
thought of as a state on $A^{**}$, and hence 
$p(\varphi) \geq 0$.  Thus $p$ gives a
nonnegative scalar function on $S(A)$, or on 
the quasistate space (that is, $\{ \alpha \varphi : 
0 \leq \alpha \leq 1 , \varphi \in S(A) \}$).
  We shall see that this function is
lowersemicontinuous if and only if $p$ is open in $A^{**}$.

In the following generalization of a well known result from
the $C^*$-algebra theory \cite{Ped}, we assume for simplicity 
that $A$ is  unital.
    If $A$ is only approximately unital then a similar result holds
with a similar proof, but one must use the quasistate space in
place of $S(A)$; this is weak* compact.  

\begin{theorem} \label{her11}  Suppose that
 $A$ is a unital-subalgebra of a $C^*$-algebra
$B$.
If $p$ is a projection in $A^{\perp \perp} \cong A^{**}$,
then the following are equivalent:
\begin{itemize}
\item [(i)]  $p$ is open as a projection in $B^{**}$ (or,
equivalently, in $A^{**}$).
\item [(ii)]  
  The set $F_p = \{ \varphi \in S(A) : \varphi(p) = 0 \}$
is a weak* closed face in $S(A)$.
 \item [(iii)]  $p$ is lowersemicontinuous on $S(A)$.
\end{itemize}
\end{theorem}

\begin{proof}
(i) $\Rightarrow$ (ii) \ 
For any projection $p \in A^{\perp \perp}$, the set 
$F_p$ is a face in $S(A)$.  For if $\psi_i \in S(A)$, $t \in [0,1]$,
and $t \psi_1 + (1-t) \psi_2 \in F_p$ then $t \psi_1(p) +
(1-t) \psi_2(p) = 0$ which forces $\psi_1(p) = \psi_2(p) = 0$,
 and $\psi_i \in F_p$.  
If $p$ is open then $G_p = \{ \varphi \in S(B): \varphi(p) = 0 \}$
is a weak* compact face in $S(B)$ by \cite[3.11.9]{Ped}.
The restriction map $r : \varphi \in S(B) \mapsto 
\varphi_{|A}  \in S(A)$ is weak* continuous, 
and maps $G_p$ into $F_p$.  On the other hand, if $\varphi \in F_p$
and $\hat{\varphi}$ is a Hahn-Banach extension of $\varphi$ to $B$
then one can show that $\langle p , \varphi \rangle 
= \langle p , \hat{\varphi} \rangle$, and so the map $r$ above 
maps $G_p$ onto $F_p$.   Hence $F_p$ is weak* closed.

(ii)  $\Rightarrow$ (i)  \ We use the notation of the
last paragraph.   If $F_p$ is weak* closed, then 
the inverse image of $F_p$ under $r$ is weak* closed.
But this inverse image is $G_p$, since 
if $\varphi \in S(B)$ then $\langle p , \varphi \rangle 
= \langle p , r(\varphi) \rangle$ by a fact in the       
last paragraph.  Thus by \cite[3.11.9]{Ped} we have (i).

 (i)  $\Rightarrow$ (iii) \ If $p$ is open,
then $p$ is a lowersemicontinuous function on $S(B)$.
Thus $\{ \varphi \in S(B) : \langle p ,  \varphi 
\rangle \leq t \}$ is weak* compact for any $t \geq 0$.
Hence its image under the map $r$ above, is 
weak* closed in $S(A)$.  However, as in the above, this
image is $\{ \varphi \in S(A) : \langle p ,  \varphi \rangle \leq t \}$.
  Thus $p$ is lowersemicontinuous  on $S(A)$. 

(iii)  $\Rightarrow$ (i) \ 
If $p$ gives a lowersemicontinuous function on $S(A)$, 
then the composition of this function with $r : 
S(B) \to S(A)$ is lowersemicontinuous on $S(B)$.
By facts in \cite[p.\ 77]{Ped}, we have that $p$ is open. 
     \end{proof}

{\bf Remark.}
  Not all weak* closed faces of $S(A)$ are 
of the form in (ii) above.  For example, let $A$ be the 
algebra of $2 \times 2$ upper triangular matrices with 
constant diagonal entries.   In this case $S(A)$ may be 
parametrized by complex numbers $z$ in a closed disk of
a certain radius centered at the origin.  Indeed states
are determined precisely by the assignment $e_{12} \mapsto z$.
The faces of $S(A)$ thus include the faces corresponding to
singleton sets of points on the 
boundary circle; and none of these faces equal $F_p$
for a projection $p \in A = A^{\perp \perp}$.        

\medskip

In view of the classical situation,
it is natural to ask about the relation between 
minimal closed projections in $B^{**}$ which
lie in $A^{\perp \perp}$ and 
the noncommutative Shilov boundary mentioned in
the introduction.   By the universal property 
of the latter object,   
if $B$ is generated as a
$C^*$-algebra by its subalgebra $A$, then there is 
a canonical
$*$-epimorphism $\theta$ from $B$ onto the 
noncommutative Shilov boundary of $A$, which in 
this case is a $C^*$-algebra.
The kernel of $\theta$ is called 
(Arveson's) {\em Shilov boundary ideal}
for $A$.  See e.g.\ \cite{SOC} and the third remark in 
\cite[4.3.2]{BLM}.
 
\begin{proposition} \label{min}  If $B$ is generated as a 
$C^*$-algebra by a closed unital-subalgebra
$A$, let $p$ be the open central
projection in $B^{**}$ corresponding to the Shilov ideal
for $A$.  Then $p^\perp$ dominates all minimal 
projections in $B^{**}$ which lie in $A^{\perp \perp}$.
\end{proposition}

\begin{proof}  Suppose that $q$ is a minimal
projection in $B^{**}$ which lies in $A^{\perp \perp}$. 
Then either $q p = 0$ or $q p = q$.
Suppose that $q p = q$.  
If $\theta$ is as above,
then since $\theta$ annihilates
the Shilov ideal we have $$\theta^{**}(q) =
\theta^{**}(q p) = \theta^{**}(q) \theta^{**}(p) = 
0 .$$
On the other hand, $\theta$ is a complete isometry from
the copy of $A$ in $B$ to the copy of $A$ in $\theta(B)$,
and so $\theta^{**}$ restricts to a complete isometry on
$A^{\perp \perp}$.   
Thus $q p = 0$, so that $q = q p^\perp$ and $q \leq p^\perp$.
\end{proof}

{\bf Example.}
The sup of closed projections in $A^{**}$ which are also minimal
projections in $B^{**}$ need not 
give the `noncommutative Shilov boundary'.
Indeed if $A$ is the $2 \times 2$ upper triangular
matrices with constant main diagonal entries, then 
there are 
no nonzero  minimal projections in $M_2$ which lie in $A$.

\section{HEREDITARY $M$-IDEALS}

A {\em left $M$-projection} of an operator space $X$ is a projection
in the $C^*$-algebra of (left) adjointable
maps on $X$; and the latter may be viewed as the restrictions
of  adjointable right module maps on a $C^*$-module
containing $X$
(see e.g.\ Theorem 4.5.15 and Section 8.4 in
\cite{BLM}).   This $C^*$-module can be taken to be
the ternary envelope of $X$.  The range of a left $M$-projection
is a {\em right $M$-summand} of $X$.
A {\em right $M$-ideal} 
of an operator space $X$ is a subspace  $J$ such that
$J^{\perp \perp}$ is a right $M$-summand of $X^{**}$.
The following result from \cite{BEZ} has been
sharpenened in the summand case:

\begin{proposition} \label{lma}  
If $A$ is an approximately unital operator algebra,
then the left $M$-projections on $A$ are precisely
`left multiplications' by projections in the multiplier algebra $M(A)$.
Such projections are all open in $A^{**}$.
The right $M$-summands of $A$ are thus the spaces $p A$
for a projection $p \in M(A)$.
The right $M$-ideals of $A$ coincide with the r-ideals of $A$.
\end{proposition}

\begin{proof}
We claim that if $p$ is a projection 
(or more generally, any hermitian) in the 
left multiplier algebra $LM(A)$,
then $p \in M(A)$.  Suppose that $B$ is a $C^*$-algebra generated by $A$,
and view $LM(A) \subset A^{\perp \perp}  \subset B^{**}$.   If $a \in A$
and if $(e_t)$ is a cai for $A$, then 
by \cite[Lemma 2.1.6]{BLM} we have $p a^* = \lim_t p  e_t a^* \in B$.
Thus $p$ is a selfadjoint element of $LM(B)$, 
and so $p \in M(B)$.  Thus $A p \subset
B \cap A^{\perp \perp} = A$, and so $p \in M(A)$.  Hence $p$ is
open as remarked early in Section 2.
The remaining assertions follow from \cite[Proposition 6.4]{BEZ}.
\end{proof}

The $M$-ideals of a unital operator algebra
are the approximately unital two-sided ideals
\cite{EfR2}.
In this case these coincide with 
the {\em complete $M$-ideals} of \cite{EfR}, which are
shown in \cite{BEZ} to be just the right $M$-ideals
which are also left $M$-ideals.   See e.g.\
\cite[Section 7]{BZ} for more information on these.
The HSA's of $C^*$-algebras are just the selfadjoint
inner ideals as remarked in the 
introduction; or equivalently as we shall see
below, they are the
selfadjoint `quasi-$M$-ideals'.
With the above facts in mind, it is
tempting to try to extend some of our results for ideals and
hereditary algebras to general $M$-ideals, be they
one-sided, two-sided, or `quasi'.   A first step along these
lines is motivated by the fact, which we have explored
in Theorem \ref{cher1} and in \cite{H}, that
r-ideals in an operator algebra
$A$ are closely tied to a matching right ideal
in a $C^*$-algebra $B$ containing $A$.   We will show
that a general (one-sided, two-sided, or `quasi')
$M$-ideal in an arbitrary operator space $X$ is
the intersection of $X$ with
the same variety of $M$-ideal in any $C^*$-algebra
or TRO containing $X$.   This generalizes a well known
fact about $M$-ideals in subspaces of $C(K)$ spaces
(see \cite[Proposition I.1.18]{HWW}).  

For an operator space $X$,
Kaneda proposed in \cite{Kun} a {\em quasi-$M$-ideal} of $X$ to be a
subspace $J \subset X$ such that $J^{\perp \perp} = p X^{**} q$
for respectively left and right $M$-projections $p$
and $q$ of $X^{**}$.   Right (resp.\ two-sided, `quasi')
$M$-ideals of a TRO or $C^*$-module are exactly the
right submodules (resp.\ subbimodules, inner ideals).
See e.g.\ \cite[p.\ 339]{BLM} and \cite{EMR,ER2}.
Here, by an inner ideal of a TRO $Z$ we mean a subspace 
$J $ with $J Z^* J \subset J$.  The assertion here 
that they coincide with the quasi $M$-ideals of $Z$ 
follows immediately from 
Edwards and R\"uttimann's characterization of 
weak* closed TRO's.  Indeed
if $J$ is an inner ideal of $Z$, then so is 
$J^{\perp \perp}$; hence \cite{ER2} gives that $J^{\perp \perp}$
is of the desired form $p Z^{**} q$.  The other direction 
follows by reversing this argument (it also may be seen as a 
trivial case of Theorem \ref{csa} below).
In fact Kaneda has considered the quasi-$M$-ideals of an
approximately unital operator algebra $A$ in 
this unpublished work \cite{Kun}.  What we will need from this is 
the following argument: If $J \subset A$ is a
quasi-$M$-ideal, then by Proposition \ref{lma} it is clear
that there exist projections $p, q \in A^{**}$
such that $J^{\perp \perp} = p A^{**} q$.  Thus
$J$ is the algebra $p A^{**} q \cap A$.

\begin{proposition}  \label{qu1}   The hereditary subalgebras of
an approximately unital operator algebra $A$ are precisely the
approximately unital quasi-$M$-ideals.
\end{proposition}

\begin{proof}  
If $J \subset A$ is a
quasi-$M$-ideal, then as we stated above,
there exist projections $p, q \in A^{**}$ 
such that $J^{\perp \perp} = p A^{**} q$, and
$J = p A^{**} q \cap A$.  If this is
approximately unital then by \cite[Proposition 2.5.8]{BLM}
$p A^{**} q$ contains a projection $e$ which is
the identity of $p A^{**} q$.  Since $e = p e q$ we have
$e \leq p$ and $e \leq q$.  So $p A^{**} q = e p A^{**} q e
= e A^{**} e$.  Thus $J = e A^{**} e \cap A$, which is
a HSA.  Conversely, if $D$ is a HSA  then 
$J^{\perp \perp} = p A^{**} p$, and so 
$J$ is a quasi-$M$-ideal.  
\end{proof}           

If ${\mathcal S}$ is a subset of a TRO we write 
$\langle {\mathcal S} \rangle$ for the subTRO
generated by ${\mathcal S}$.  We write $\widehat{ \phantom{.} } \,$ 
for the canonical map from a space into its second dual.

\begin{lemma} \label{stros}  If $X$ is an operator space, and
 if $({\mathcal T}(X^{**}),j)$  is a ternary envelope of
$X^{**}$, then $\langle j(\hat{X}) \rangle$ is a ternary envelope of $X$.
\end{lemma}

\begin{proof}
This follows from a diagram chase.  Suppose that $i : X \to W$ is
a complete isometry into a TRO, such that $\langle i(X) \rangle = W$.  
Then $i^{**} : X^{**} \to W^{**}$ is
a complete isometry.  By the universal property of the 
ternary envelope, there is a ternary morphism 
$\theta : \langle i^{**}(X^{**}) \rangle \to {\mathcal T}(X^{**})$
 such that $\theta \circ i^{**} = j$.  Now $W$ may also be regarded 
as the subTRO of $W^{**}$ generated by $i(X)$, 
and the restriction $\pi$ of
$\theta$ to $W = \langle i(X) \rangle$ is a 
ternary morphism into ${\mathcal T}(X^{**})$ which has the 
property that $\pi(i(x)) = j(\hat{x})$.  Thus
$\langle j(\hat{X}) \rangle$ has the universal property of the
ternary envelope.  
 \end{proof}

\begin{theorem}  \label{csa}  Suppose that $X$ is a subspace of a TRO $Z$
and that $J$ is a right $M$-ideal (resp.\ quasi $M$-ideal,
complete $M$-ideal) of $X$.
In the `complete $M$-ideal' case we also assume that $\langle X \rangle
= Z$.
Then $J$ is the intersection of $X$
with the right $M$-ideal (resp.\ quasi $M$-ideal,
complete $M$-ideal) $J Z^* Z$ (resp.\ $J Z^* J$,
$Z J^* Z$) of $Z$.  \end{theorem}

\begin{proof}
There are three steps.  We will also use the fact
that in a TRO $Z$, for any $z \in Z$ we have that $z$ lies in the
closure of $z \langle z \rangle^* z$.
 This follows by considering the polar decomposition 
$z = u |z|$, which implies that 
$z z^* z = u |z|^3$, for example.  Then use 
the functional calculus for $|z|$, and the 
fact that one may approximate the monomial
$t$ by polynomials in $t$ with only odd powers and degree $\geq 3$.  
Similarly, $z$ lies in the
closure of $\langle z \rangle z^* \langle z \rangle$.

First, suppose that $Z$
is the ternary envelope ${\mathcal T}(X)$ of $X$.
Suppose that $J$ is a right $M$-ideal (resp.\
quasi-$M$-ideal, complete  $M$-ideal)
  in  $X$. 
 If $({\mathcal T}(X^{**}),j)$  is a ternary envelope of
$X^{**}$, then $j(J^{\perp \perp}) = p j(X^{**})$ for
a left adjointable projection $p$ (resp.\ $j(J^{\perp \perp}) = p j(X^{**}) q$
for left/right adjointable projections $p, q$) on  ${\mathcal T}(X^{**})$.
In the complete  $M$-ideal case we have $p w = w q$ for all 
$w \in {\mathcal T}(X^{**})$; this follows from e.g.\ 
\cite[Theorem 7.4 (vi)]{BZ} and its proof.
  We view ${\mathcal T}(X) \subset {\mathcal T}(X^{**})$ as above.
Let $\tilde{J}$ be the set of
$z \in {\mathcal T}(X)$ such that $p z = z$ (resp.\ $z = p z q$).
 Then $\tilde{J} \cap j(\hat{X}) = j(\hat{J})$,
since $J = J^{\perp \perp} \cap X$.
Next, define $\bar{J} = j(\hat{J}) {\mathcal T}(X)^* {\mathcal T}(X)$
(resp.\ $= j(\hat{J}) {\mathcal T}(X)^* j(\hat{J})$,
 $=   {\mathcal T}(X)  j(\hat{J})^* {\mathcal T}(X))$.   This
is a right $M$-ideal (resp.\ inner ideal,
$M$-ideal)  in ${\mathcal T}(X)$,
and it is clear, 
using the fact in the first paragraph of the proof, that
$$j(\hat{J}) \subset \bar{J} \cap j(\hat{X}) \subset \tilde{J}
\cap j(\hat{X}) = j(\hat{J}) .$$
Thus $J = \bar{J} \cap X$.

In the rest of the proof we consider
only the quasi $M$-ideal case, the others are similar.

Second, suppose that $X$ generates $Z$ as a TRO.
Let $j : X \to {\mathcal T}(X)$ be the Shilov embedding.
If $x \in (J Z^* J)  \cap X$ then  applying the
universal property of ${\mathcal T}(X)$ 
there exists a ternary morphism $\theta : Z \to {\mathcal T}(X)$
with 
$$j(x) = \theta(x) \in j(X) \cap \theta(J) \theta(Z)^* \theta(J) \subset
j(X) \cap j(J) {\mathcal T}(X)^* j(J) =
j(X) \cap \bar{J} = j(J) ,$$ by the
last paragraph.  Hence $x \in J$.

Third, suppose that $X \subset Z$, and that the subTRO
$W$ generated by $X$ in $Z$ is not $Z$.  We claim that
$X \cap (J Z^* J) = X \cap (J W^* J)$.
To see this, we set $J' = J W^* J$.  This
is an inner ideal in $W$.  Moreover,  $J \subset J'$ by the
fact at the start of the proof.
We claim that for any inner ideal $K$ in $W$,
we have $(K Z^* K) \cap W = K$.  Indeed if $e$ and $f$ are
the support projections for $K$, then
$$(K Z^* K) \cap W \subset (e Z f) \cap W \subset (e W f) \cap W = K .$$
This implies
$$X \cap (J Z^* J) \subset X \cap (J' Z^* J')
\subset X \cap J' = X \cap (J  W^* J) = J,$$
as required.
\end{proof}

\section{REMARKS ON PEAK AND $p$-PROJECTIONS}

Let $A$ be a unital-subalgebra of a $C^*$-algebra $B$.
We recall from \cite{H} that a {\em peak projection} $q$ for $A$ is
a closed projection in $B^{**}$, such that there exists
an $a \in {\rm Ball}(A)$ with  $q a = q$ and 
 satisfying any one of a long list 
of equivalent conditions; for example $\Vert a r \Vert < 1$ for 
every closed projection $r \leq q^\perp$.   
  We say that $a$ peaks at $q$.  A {\em  $p$-projection}
 is an infimum of peak projections;
 and this is equivalent to it being a 
weak* limit of a decreasing net of peak projections
by \cite[Proposition 5.6]{H}.   Every 
$p$-projection is an {\em approximate $p$-projection},
where the latter term means
a closed projection in $A^{**}$.   The most glaring problem 
concerning these projections is that it is currently unknown whether
the converse of this is true, as is the
case in the classical setting of function algebras
\cite{Gam}.   Motivated partly by this question,
 in this section we offer several 
results concerning these projections. 
   Our next result implies that
this question is equivalent to the following simple-sounding question:

\smallskip

{\em Question:} Does every approximately unital operator algebra $A$
have an approximate identity
of form $(1- x_t)$ with  $x_t \in {\rm Ball}(A^1)$? 
Here $1$ is the identity of the unitization 
$A^1$ of $A$.

\medskip

Equivalently, does every operator algebra $A$
with a left cai have a left cai of the form $(1- x_t)$ for $x_t \in {\rm Ball}(A^1)$?

\medskip

By a routine argument,
these are also equivalent to: If $A$ is an approximately unital operator algebra
and $a_1, \cdots , a_n \in A$ and 
$\epsilon > 0$, does there exist $x \in {\rm Ball}(A^1)$ with
$1-x \in A$ and $\Vert 
x  a_k \Vert < \epsilon$ for all $k = 1, \cdots, n$? 

\medskip

Note that if these were true, and if $A$ does not have an identity,
then necessarily $\Vert x_t \Vert = 1$. For if $\Vert x_t \Vert < 1$ 
then $1- x_t$ is invertible in $A^1$, so that
$1 \in A$.

\begin{theorem} \label{one} If $J$ is a closed subspace of
a unital operator algebra  $A$, then
the following are equivalent:
\begin{itemize}
\item [(i)]  $J$ is a right ideal with a left approximate identity 
(resp.\  a HSA with approximate identity) of  
the form $(1- x_t)$ for $x_t \in {\rm Ball}(A)$.
\item [(ii)] $J$ is an r-ideal 
(resp.\ HSA) for whose support projection $p$ we
have  
that $p^\perp$ is a $p$-projection for $A$.
  \end{itemize}
 \end{theorem}

\begin{proof} 
 Suppose that $J = \{ a \in A : q^\perp a = a \}$ for a $p$-projection $q$
for $A$ in $B^{**}$.  We may suppose that $q$ is a decreasing weak* limit
of a net of peak projections $(q_t)$ for $A$.  
If $a \in A$ peaks at $q_t$, then
by a result in \cite{H} we have that $a^n \to q_t$ weak*.
Next let ${\mathcal C} = \{ 1 - x : x \in {\rm Ball}(A) \} \cap J$,
a convex subset of $J$ containing the elements $1-a^n$ 
above.  Thus $q_t^\perp \in \overline{{\mathcal C}}^{w*}$
and therefore $q^\perp \in \overline{{\mathcal C}}^{w*}$.
Let $e_t \in {\mathcal C}$ with $e_t \to q^\perp$ w*.
Then $e_t x \to q^\perp x = x$ weak* for all $x \in J$.
Thus $e_t x \to x$ weakly.   Next, for fixed
$x_1, \cdots , x_m \in J$ consider the convex set
$F = \{ (x_1 - u x_1 , x_2  - u x_2 , \cdots , x_m  - u x_m ) :
u \in {\mathcal C} \}$.   (In the HSA case 
one also has to include coordinates $x_k - x_k u$ here.) 
Since $(0,0,\cdots ,0)$ is in the
weak closure of $F$ it is in the norm closure.  
Given $\epsilon > 0$, there exists $u \in {\mathcal C}$ such that
$\Vert x_k - u x_k \Vert < \epsilon $ for all $k = 1, \cdots , m$.
From this it is clear (see the end of the proof of \cite[Proposition 2.5.8]{BLM})
that there is a left approximate identity
for $J$ in ${\mathcal C}$, which  shows (i).

Suppose that $J$ is a right ideal with a left approximate
identity $(e_t)$ of the stated form $e_t = 1 - x_t$.
If $(x_{t_\mu})$ is any w*-convergent subnet of
$(x_t)$, with limit $r$, then $\Vert r \Vert \leq 1$.
Also $1 - x_{t_\mu} \to 1 - r$.   On the other hand,
$(1 - x_{t_\mu})x \to x$ for any $x \in J$, so that
$(1-r) x = x$.  Hence $(1-r) \eta = \eta$ for any $\eta \in J^{\perp \perp}$,
so that $1-r$ is the (unique) left identity $p$ for $J^{\perp \perp}$.
Hence $1-r$ is idempotent, so that $r$ is idempotent.  Hence
$r$ is an orthogonal projection, and therefore so also is
$p = 1-r$.   Also, $e_t \to p$ w*, by a fact in topology about
nets with unique accumulation points.
We have $J = p A^{**} \cap A =
\{ a \in A : p a = a \}$.    Since $p$ has norm $1$,
$J$ has a left cai.    Since $p e_t = e_t$,
$p$ is an open projection in $B^{**}$, so that $q = 1-p$ is closed.
If $a = e_t$ then we have that $l(a)^\perp (1-a) = l(a)^\perp$,
where $l(a)$ is the left support projection for $a$.
Thus by a result in \cite{H}  there is a peak projection $q_a$ with peak
 $a_0 =  1 - \frac{a}{2} \in A$
such that $l(a)^\perp \leq q_a$.   Since $a_0^n \to q_a$ weak*,
and since $(1-p) a_0^n = 1-p$, we have $(1-p) q_a = 1-p$.  That is,
$q \leq q_a$.   Let $J_a = \{ x \in A :  q_a^\perp x = x \}$.
By the last paragraph, $J_a$ is an r-ideal,
and since $q \leq q_a$ we have that $J_a \subset J$.
The closed span of all the $J_a$ for $a = e_t$
equals $J$, since $e_t \in J_{e_t}$ and
any $x \in J$ is a limit of $e_t x \in J_{e_t}$.    
By the proof of
\cite[Theorem 4.8.6]{BLM} we deduce that the supremum of the $q_a^\perp$
equals $q^\perp$.   Thus $q$ is a $p$-projection.
The HSA case follows easily from this and Corollary \ref{her7}.
  \end{proof}

\begin{corollary}  \label{two} 
Let $A$ be a unital-subalgebra of a $C^*$-algebra $B$.
A projection $q \in B^{**}$ is a $p$-projection for
$A$ in $B^{**}$, if and only if there exists a net $(x_t)$ in ${\rm Ball}(A)$ with
$q x_t = q$, and $x_t \to q$ weak*.
\end{corollary}

\begin{proof}  Supposing that $q$ is a $p$-projection,
we have by the last result
that  $J = \{ a \in A : q^\perp a = a \}$  has a left approximate
identity $(1 - x_t)$ with $x_t \in {\rm Ball}(A)$,
and by the proof of that result $q^\perp$ is the
support projection, so that $1 - x_t \to q^\perp$ weak*.

Conversely,  supposing the existence of such a
net, let $J = \{ a \in A : q^\perp a = a \}$.   This is a
right ideal.  Moreover $J^{\perp \perp} \subset q^\perp A^{**}$.
 If $a \in A$ then $q^\perp a = \lim_t \, (1 - x_t) a
\in J^{\perp \perp}$.  By a similar argument,
$q^\perp \eta \in J^{\perp \perp}$ for any $\eta \in A^{**}$.
Thus  $J^{\perp \perp} = q^\perp A^{**}$, and so
$q^\perp$ is the support projection for $J$, and $J$ has a
left cai.   By a slight variation of
the argument at the end of the first
paragraph of the proof of the last result,
$J$ satisfies (i) of that result, and hence
by that result $q$ is a $p$-projection.
\end{proof}

The following known result (see e.g.\ \cite{MAII,DP}) is quite 
interesting in light of the 
question just above Theorem \ref{one}.

\begin{proposition} \label{nine}  If $J$ is an
nonunital operator algebra with a cai (resp.\ left cai), 
then $J$ has an (resp.\ a left) approximate identity
of the form  $(1-x_t)$, where $x_t \in J^1$ and
$\lim_t \, \Vert x_t \Vert = 1$ and 
$\lim_t \, \Vert 1-x_t \Vert = 1$.
Here $J^1$ is the unitization of $J$.
\end{proposition}

\begin{proof}   We just sketch the proof in the 
left cai case, following the proof of \cite[Theorem 3.1]{MAII}. 
  Let $A = J^1$.  Thus $J$ is an
r-ideal in the unital operator algebra $A$.
Suppose that the support projection is
$p = q^\perp \in A^{**}$, and that $(u_t)$ is the left cai in 
$J$.  If $B$ is a $C^*$-algebra generated by $A$, then there is
an increasing net in ${\rm Ball}(B)$ with 
weak* limit $p$.   
We can assume that the increasing net is indexed by the same directed set.
Call it $(e_t)$.  Since $e_t - u_t \to 0$ weakly, new nets of convex combinations
$(\widetilde{e_s})$ and $(\widetilde{u_s})$ will satisfy 
$\Vert \widetilde{e_s} - \widetilde{u_s} \Vert \to 0$.
We can assume that $(\widetilde{u_s})$ is a left cai for $J$.
We have 
$$\Vert 1 - \widetilde{u_s} \Vert  \leq \Vert 1 - \widetilde{e_s} \Vert
+ \Vert  \widetilde{e_s} - \widetilde{u_s} \Vert  \leq 1 +
\Vert  \widetilde{e_s} - \widetilde{u_s} \Vert  \to 1 .$$ 
The result follows easily from this.  \end{proof}

We are also able to give another  
characterization of $p$-projections,
which is of `nonselfadjoint Urysohn lemma' or `peaking' 
flavor, and therefore should be useful in  
future applications of
`nonselfadjoint peak interpolation'.    This result should be compared 
with \cite[Theorem 5.12]{H}.

\begin{theorem} \label{peakch}  Let $A$ be a
unital-subalgebra of $C^*$-algebra $B$ and let
$q \in B^{**}$ be a closed projection.  Then 
$q$ is a $p$-projection for $A$ iff for any open
projection $u \geq q$, and any $\epsilon > 0$,
there exists an $a \in {\rm Ball}(A)$ with 
$a q = q$ and $\Vert a (1-u) \Vert < \epsilon$
and $\Vert (1-u)  a \Vert < \epsilon$.
\end{theorem}

 \begin{proof}
($\Leftarrow$) \ This follows by an easier variant of
 the proof of  \cite[Theorem 4.1]{H}.
Suppose that for each open $u \geq q$, and
positive integer $n$, there exists an $a_n \in {\rm Ball}(A)$ with
$a_n q = q$ and $\Vert a_n (1-u) \Vert < 1/n$.   
By taking a weak* limit we find $a \in A^{\perp\perp}$
with $a q = q$ and $a (1-u) = 0$. 
We continue  as in \cite[Theorem 4.1]{H}.
Later in the proof where $q_n$ is defined, we appeal
to Lemma 3.5 in place of Lemma 3.6, so that $q_n$ is a
peak projection.  Now the proof is quickly finished:
Let $Q = \bigwedge_n q_n$, a $p$-projection.  As in
the other proof we have that $q \leq Q \leq r \leq u$, and
that this forces $q = Q$.  Thus $q$ is a $p$-projection.

($\Rightarrow$) \ Suppose that $q$ is a $p$-projection,
and $u \geq q$ with $u$ open.  By `compactness' of $q$
 (see the remark just above \cite[Proposition 2.2]{H}), 
there is a peak projection $q_1$ with 
$q \leq q_1 \leq u$.  Note that if $a q_1 = q_1$ then
$a q = a q_1 q = q_1 q = q$.   Thus we may assume that 
$q$ is a peak projection.  By the noncommutative Urysohn lemma
\cite{Ake2}, there is an $x \in B$ with 
$q \leq x \leq u$.   Suppose that  $a \in {\rm Ball}(A)$
peaks at $q$, and $a^n \to q$ weak* (see e.g.\ \cite[Lemma 3.4]{H}
or the results below). 
Then  $a^n (1-x) \to q (1-x) = 0$ weak*, and hence weakly
in $B$.  Similarly, $(1-x) a^n \to  0$ weakly.
By a routine convexity argument in $B \oplus B$, given $\epsilon > 0$
there is a convex combination $b$ of the $a^n$
such that $\Vert b (1-x) \Vert < \epsilon$ and 
$\Vert  (1-x) b \Vert < \epsilon$.
Therefore $\Vert b (1-u) \Vert = 
\Vert b (1-x) (1-u)  \Vert < \epsilon$.
Similarly $\Vert   (1-u) b \Vert  < \epsilon$. \end{proof}

We would guess that being a $p$-projection is also equivalent to 
the special case where $a = 1$ and $x \leq 1$ of the following
definition. 

If $A$ is a unital-subalgebra of $C^*$-algebra $B$
and if $q \in B^{**}$ is closed then we say that $q$ is
a {\em strict $p$-projection} if given $a \in A$ and a
strictly positive $x \in B$ with $a^*qa \le x$, then
there exists $b
\in A$ such that $qb=qa$ and $b^*b \le x$.   In \cite[Proposition 3.2]{H}
it is shown that if $q$ is a projection
in $A^{\perp\perp}$ then
the conditions in the last line hold  except that
$b^*b \le x + \epsilon$.
So being a strict $p$-projection is the case $\epsilon = 0$ of
that interpolation result.

\begin{corollary} \label{strict}  Let $A$ be a
unital-subalgebra of $C^*$-algebra $B$ and let
$q \in B^{**}$ be a strict $p$-projection for $A$.
Then $q$ is a $p$-projection.
\end{corollary}
 
\begin{proof}
Using the noncommutative Urysohn lemma
as in the first 
 few lines of the proof of  \cite[Theorem 4.1]{H},
it is easy to see that $q$ satisfies the condition in 
Theorem \ref{peakch}.
\end{proof}

The above is related to
the question of whether every r-ideal $J$ in
a (unital say) operator algebra is `proximinal'
(that is, whether every $x \in A$ has a closest point in $J$). 
 
\begin{proposition} \label{prox}   If
$q$ is a strict $p$-projection for
a unital operator algebra $A$,
 then the corresponding r-ideal $J =
q^\perp A^{**} \cap A$ is proximinal in $A$.
\end{proposition}
 
\begin{proof}
Let $a \in A$.  By Proposition 3.1 in \cite{H}, $\|a + J\| =
\|qa\|$.  Also, $a^*qa \le \|qa\|^2$, so by hypothesis there exists $b
\in A$ such that $qb = qa$ and $b^*b \le \|qa\|^2$.  Thus $\|b\|^2 =
\|b^*b\| \le \|q a\|^2$.
Then $\|a + J\| = \|qa\|
\geq \|b\| = \|a + (b-a)\|$.  However,
$b-a \in J$ since $q(b-a) = 0$.  So $J$ is proximinal.
\end{proof}

Some of the results below stated for right ideals also have
HSA variants which we leave to the reader.

\begin{proposition}  \label{four} A $p$-projection $q$ for 
a unital operator algebra $A$ is a
peak projection iff the associated right ideal is of the
form $\overline{(1-a) A}$ for some $a \in {\rm Ball}(A)$.
In this case, $q$ is the peak for $(a+1)/2$.    \end{proposition}

\begin{proof}
Let $J =  \{ a \in A : q^\perp a = a \}$, for a $p$-projection $q$.

($\Rightarrow$) \ If $q$ peaks at $a$ then $q^\perp (1-a) = (1-a)$, so that
$(1-a)A \subset J$.  If $\varphi \in ((1-a)A)^\perp$ then
$\varphi((1-a^{n+1}) A) = \varphi((1-a)(1 + a + \cdots + a^n) A) = 0$.
In the limit we see that $\varphi \in (q^\perp A)_\perp$, so that
$\varphi \in J^\perp$.  Hence $J = \overline{(1-a) A}$.

($\Leftarrow$) \ Suppose that $J = \overline{(1-a) A}$ for some $a \in {\rm Ball}(A)$.
Then $q a = q$, and by a result in \cite{H} there exists
a peak projection $r \geq q$ with peak $b = (a+1)/2$.
Since $1 - b = (1-a)/2$, it is clear that $J = \overline{(1-b) A}$.
If $(e_t)$ is the left cai for $J$ then  $r^\perp e_t = e_t$.
In the limit, $r^\perp q^\perp = q^\perp$, so that $r \leq q$.
Thus $r = q$.  \end{proof}

This class of `singly generated' right ideals has played an important
role in some work of G. A. Willis (see e.g.\  \cite{Wil}).

\begin{lemma} \label{wilp}  If $A$ is an operator algebra, 
and if $a \in {\rm Ball}(A)$
then $\overline{(1-a) A}$ is an r-ideal of $A$
with a sequential left approximate identity
of the form $(1- x_n)$ for $x_n \in {\rm Ball}(A)$.
Similarly, $\overline{(1-a) A (1-a)}$ is a HSA of $A$.
\end{lemma}

\begin{proof}  Let $J = \overline{(1-a) A}$,
and let $e_n =  1 - \frac{1}{n} \sum_{k=1}^n a^k$,
which is easy to see is in $(1-a) A$.
Moreover,
$$e_n (1-a) = 1 - \frac{1}{n} \sum_{k=1}^n a^k
- a + \frac{1}{n} \sum_{k=2}^{n+1} a^k =
1 - a - \frac{1}{n}(a - a^{n+1}) \to 1-a .$$
Note that $J$ is an r-ideal 
by Theorem \ref{one}.  We leave the rest to the reader.    \end{proof}

\begin{corollary}  \label{seven}  If
$a$ is a contraction in a unital $C^*$-algebra $B$ then
\begin{itemize}
\item [(i)]  The Cesaro averages of $a^n$ converge weak* to a
peak projection $q$ with $q a = q$.
\item [(ii)]  If $a^n \to q$ weak* then $q$ is a  peak projection.
Conversely, if $q$ is a  peak projection then there 
exists an $a \in {\rm Ball}(B)$ with $a^n \to q$ weak*.
\end{itemize}
 Also $q$ is the peak for $(a+1)/2$.
  \end{corollary}

\begin{proof}   (i) \ By Theorem \ref{one} and Lemma \ref{wilp}
(and its proof),
$J = \overline{(1-a) A} = \{ a \in A : q^\perp a = a \}$
for a $p$-projection $q$ which is a weak* limit
of $e_n = 1 - \frac{1}{n} \sum_{k=1}^n a^k$.
Thus $\frac{1}{n} \sum_{k=1}^n a^k \to q$ weak*, and
clearly $q a = q$.  By \ref{four} and its proof, $q$ is a peak projection
with $(a+1)/2$ as a peak.

(ii) \ If $a^n \to q$ weak* then it is easy to check that
$\frac{1}{n} \sum_{k=1}^n a^k  \to q$ weak*.  Thus one direction
of (ii) follows from (i), and the other direction
is in \cite{H}.  \end{proof}

{\bf Remarks.}
1)  \  In fact it is not hard to show that
the Cesaro averages in (i) above converge 
strongly, if $B$ is in its universal representation.   

\medskip 

2)  \ We make some remarks on support projections.
We recall from \cite{H} that
if $q$ is a projection in
$B^{**}$ and if $q$ peaks at a contraction $b \in B$
then $q^\perp$ is the right support projection $r(1-b)$.
Conversely, if $b \in {\rm Ball}(B)$ then
the complement of the right support projection $r(1-b)$
is a peak projection which peaks at $(1+b)/2$.
Thus the peak projections are precisely the
complements of the
right support projections $r(1-b)$ for contractions $b \in B$.

It follows that $q$ is a $p$-projection for a
unital-subspace $A$ of a $C^*$-algebra
$B$ iff $q = \wedge_{x \in {\mathcal S}} \,
r(1-x)^\perp$ for a nonempty subset ${\mathcal S} \subset
{\rm Ball}(A)$.

Also, if $J$ is a right ideal of 
a unital operator algebra $A$, and if $J$ has a left approximate identity
of the form $(1 - x_t)$ with $x_t \in {\rm Ball}(A)$,
then it is easy to see that the support projection of $J$ is 
$\vee_t \, l(1 - x_t)$.


\begin{thebibliography}{99}
\bibitem{Ake}
\textsc{C.A. Akemann}, The general Stone-Weierstrass problem,
\textit{J. Funct. Anal.} \textbf{4}(1969), 277--294.

\bibitem{Ake2} \textsc{C.A. Akemann}, Left ideal structure of
    $C^*$-algebras, \textit{J. Funct. Anal.} \textbf{6}(1970), 305--317.

\bibitem{FSIO} \textsc{C.A. Akemann, G.K. Pedersen}, Facial
structure in operator algebra theory, \textit{Proc. London Math.
Soc.} \textbf{64}(1992), 418--448.

\bibitem{MAII} \textsc{A. Arias, H.P. Rosenthal}, $M$-complete
approximate identities in operator spaces, \textit{Studia Math.}
 \textbf{141}(2000), 143--200.

\bibitem{SOC}  \textsc{W.B. Arveson}, Subalgebras of $C^{*}$-algebras,
\textit{Acta Math.} \textbf{123}(1969), 141--224.

\bibitem{Bghm}  \textsc{D.P. Blecher},  A generalization of Hilbert
modules, \textit{J. Funct. Anal.} \textbf{136}(1996), 365--421.

\bibitem{Bnat}  \textsc{D.P. Blecher}, Some
general theory of operator algebras and their modules, in
\textit{ Operator Algebras and Applications}, Ed. A. Katavolos,
Nato ASI Series, Series C - Vol. 495, Kluwer, Dordrecht-Boston-London 1997.

\bibitem{Bna} \textsc{D.P. Blecher}, A new approach to $C^*$-modules,
  \textit{Math. Annalen} \textbf{307}(1997), 253--290.

\bibitem{B}  \textsc{D.P. Blecher}, One-sided ideals and approximate
  identities in operator algebras,
\textit{J. Australian Math. Soc.} \textbf{76}(2004), 425--447.

\bibitem{BEZ}  \textsc{D.P. Blecher, E.G. Effros, V. Zarikian},
One-sided $M$-ideals and multipliers in operator spaces, I,
\textit{Pacific J. Math.} \textbf{206}(2002), 287--319.

\bibitem{BLM}  \textsc{D.P. Blecher, C. Le Merdy}, \textit{ Operator
algebras and their modules---an operator space approach,} Oxford Univ.\  
Press, Oxford 2004.

\bibitem{BMP} \textsc{D.P. Blecher, P.S. Muhly, V.I. Paulsen},
Categories of operator modules (Morita equivalence and projective
modules),  \textit{Mem. Amer. Math. Soc.} \textbf{143}(2000), no.\ 681. 


\bibitem{BZ} \textsc{D.P. Blecher, V. Zarikian},
The calculus of one-sided $M$-ideals and multipliers  in
  operator spaces, \textit{Mem. Amer. Math. Soc.} \textbf{179}(2006),
no.\ 842.

\bibitem{DP} \textsc{K.R. Davidson, S.C. Power}, Best
approximation in $C^*$-algebras, \textit{J. Reine Angew. Math.} 
\textbf{368}(1986), 43--62.

\bibitem{EMR}  \textsc{C.M. Edwards, K. McCrimmon, G.T. R\"uttimann},
  The range of a structural projection,
\textit{J. Funct. Anal.}  \textbf{139}(1996), 196--224.

\bibitem{ER2}  \textsc{C.M. Edwards, G.T. R\"uttimann},
On inner ideals in ternary algebras, \textit{Math. Z.}
 \textbf{204}(1990), 309--318.

\bibitem{ER1} \textsc{C.M. Edwards, G.T. R\"uttimann}, Inner ideals in
  C*-algebras, \textit{Math. Ann.} \textbf{290}(1991), 621--628.

\bibitem{ER3}  \textsc{C.M. Edwards, G.T. R\"uttimann},
A characterization of inner ideals in JB*-triples, \textit{
  Proc. Amer. Math. Soc.} \textbf{116}(1992), 1049--1057.

\bibitem{Ef} \textsc{E.G. Effros}, Order ideals in a $C^*$-algebra and
  its dual, \textit{Duke Math. J.} \textbf{30}(1963), 391--412.

\bibitem{EfR2} \textsc{E.G. Effros, Z-J. Ruan},
On non-self-adjoint operator algebras, \textit{Proc. Amer. Math. Soc.} 
\textbf{110}(1990), 915--922.

\bibitem{EfR} \textsc{E.G. Effros, Z-J. Ruan},
Mapping spaces and liftings for operator spaces, \textit{Proc. London
  Math. Soc.} \textbf{69}(1994), 171--197.

\bibitem{FR}  \textsc{Y. Friedman, B. Russo},
 Structure of the predual of a $JBW\sp *$-triple,
\textit{J. Reine Angew. Math.} \textbf{356}(1985), 67--89.

\bibitem{Gam}  \textsc{T.W. Gamelin}, \textit{Uniform algebras},
  Prentice-Hall, New York 1969.

\bibitem{Ham}  \textsc{M. Hamana}, Triple envelopes and Silov boundaries
  of operator spaces, \textit{Math. J. Toyama University}
 \textbf{22}(1999), 77--93.

\bibitem{HWW} \textsc{P. Harmand, D. Werner, W. Werner},
\textit{$M$-ideals in Banach spaces and Banach algebras},
Lecture Notes in Math., vol.\ 1547, Springer-Verlag, Berlin--New York 1993.

\bibitem{HTh}  \textsc{D.M. Hay}, \textit{Noncommutative topology and
    operator algebras} (tentative title), PhD Dissertation, University
  of Houston, 2006.

\bibitem{H}  \textsc{D.M. Hay}, Closed projections and peak interpolation
for operator algebras, Preprint (2005), math.OA/0512353.

\bibitem{Kth} \textsc{M. Kaneda}, \textit{ Multipliers and algebrizations
of operator spaces,} PhD Dissertation, University of Houston, 2003.

\bibitem{Kun} \textsc{M. Kaneda}, Unpublished note on quasi-$M$-ideals
  (2003).

\bibitem{K} \textsc{M. Kaneda}, Extreme points of the unit ball of a
  quasi-multiplier space, Preprint (2004),  Math.OA/0408235.

\bibitem{Kat} \textsc{E.G. Katsoulis}, Geometry of the unit ball and
  representation theory for operator algebras, \textit{Pacific
    J. Math.} \textbf{216}(2004), 267--292.

\bibitem{MZ}  \textsc{J-P. Ma, L-C. Zhang}, A characterization of inner
  ideals of $L^\infty$-Banach algebras and a theorem on isomorphisms between
them (Chinese), \textit{Nanjing Daxue Xuebao Ziran Kexue Ban}
 \textbf{35}(1999), 502-504.

\bibitem{Ped} \textsc{G.K. Pedersen}, \textit{$C^*$-algebras
and their automorphism groups}, Academic Press, London 1979.

\bibitem{Wil} \textsc{G.A. Willis}, Factorization in finite codimensional
ideals of group algebras, \textit{Proc. London Math. Soc.}
 \textbf{82}(2001), 676-700.

\bibitem{ZZ}  \textsc{L-C. Zhang, X. Zhang}, The characterization of inner
  ideal in concrete operator algebra and an isomorphism theorem
 between them (Chinese), \textit{Acta Mathematica Sinica} 
\textbf{43}(2000), 843-846.

\end{thebibliography}
  \end{document}